\documentclass[11pt,reqno]{amsart}
\usepackage[margin=1in]{geometry}
\usepackage{amsmath, amssymb, amsthm, thmtools, mathtools, comment}
\usepackage{booktabs}
\usepackage{enumitem}
\usepackage{microtype}
\usepackage{bbm}
\usepackage[colorlinks=true,linkcolor=blue,citecolor=blue,urlcolor=blue]{hyperref}
\usepackage{tikz}
\usetikzlibrary{positioning}
\usepackage{ytableau}
\usepackage{cleveref}

\ytableausetup{
  smalltableaux,
  boxsize=1.2em,
  centertableaux,
  mathmode
}

\numberwithin{equation}{subsection}
\newtheorem{theorem}[equation]{Theorem}
\newtheorem{lemma}[theorem]{Lemma}
\newtheorem{proposition}[theorem]{Proposition}

\theoremstyle{definition}
\newtheorem{definition}[theorem]{Definition}
\newtheorem{example}[theorem]{Example}
\newtheorem*{acknowledgement}{Acknowledgements}
\newtheorem*{aidisclosure}{AI disclosure}

\theoremstyle{remark}
\newtheorem{remark}[theorem]{Remark}

\DeclareMathOperator{\Ass}{Ass}
\DeclareMathOperator{\MaxAss}{MaxAss}

\DeclareMathOperator{\GL}{GL}
\DeclareMathOperator{\SL}{SL}
\DeclareMathOperator{\Orth}{O}

\DeclareMathOperator{\Sym}{Sym}

\newcommand{\kfield}{\Bbbk}

\newcommand{\cR}{\mathcal R}

\newcommand{\fp}{\mathfrak p}

\newcommand{\fm}{\mathfrak m}

\newcommand{\Schur}{\mathbf S}

\renewcommand{\l}{\lambda}

\title[Invariant ideals on determinantal rings]{Symbolic Powers and Asymptotic Invariants of \\ Invariant Ideals in Determinantal Rings}

\author{Sankhaneel Bisui}
\address{Department of Mathematics\\
Indian Institute of Science Education and Research Bhopal\\
Bhopal Bypass Road, Bhauri, Bhopal 462066, Madhya Pradesh, India}
\email{sankhaneel@iiserb.ac.in}

\author{Aryaman Maithani}
\address{Department of Mathematics\\
 University of Utah\\
  155 South 1400 East, Salt Lake City, UT 84112, USA}
\email{maithani@math.utah.edu}

\author{Kesavan Mohana Sundaram}
\address{Department of Mathematics\\
University of Nebraska--Lincoln\\
203 Avery Hall, Lincoln, NE 68588, USA}
\email{km2@huskers.unl.edu}

\author{Alexandra Seceleanu}
\address{Department of Mathematics\\
University of Nebraska--Lincoln\\
203 Avery Hall, Lincoln, NE 68588, USA}
\email{aseceleanu@unl.edu}

\author{Stephen Stern}
\address{Department of Mathematics\\
University of Nebraska--Lincoln\\
203 Avery Hall, Lincoln, NE 68588, USA}
\email{sstern2@huskers.unl.edu}


\begin{document}

\begin{abstract}
    Let $A$ be the coordinate ring of a bounded-rank locus of generic,
symmetric, alternating, or Hankel matrices. Thus $A$ is the quotient of a polynomial ring by a generic determinantal, symmetric determinantal, generic Pfaffian, or Hankel determinantal ideal.  In each case there is a
nested family of rank primes
$\fp_1\supset \fp_2\supset\cdots\supset\fp_q$ and a standard-monomial
basis whose shapes are diagrams with parts at most $q$.  We isolate the
common properties of these four straightening theories and use them to
give a uniform treatment of products
$D_\sigma=\fp_{\sigma_1}\cdots\fp_{\sigma_s}$ in the quotient ring $A$, extending analogous results over suitable polynomial rings.  In all four cases the shape $\sigma$ must have parts of size at most $q$.
We determine their primary decompositions and associated primes in the
rank-$q$ quotient, compute all symbolic powers with respect to the full
set of associated primes, and prove that their symbolic Rees algebras
are finitely generated. In the generic, symmetric, and alternating cases the respective symbolic Rees algebras are also normal, and Cohen--Macaulay.  We also determine  the asymptotic resurgence
$\widehat\rho$ and Waldschmidt constants of these ideals and prove they are related by $\widehat\rho(D_\sigma)=\alpha(D_\sigma)/\widehat\alpha(D_\sigma)$. 
\end{abstract}

\maketitle
\setcounter{tocdepth}{1}
\tableofcontents

\section{Introduction}\label{sec:introduction}

Determinantal varieties exhibit a natural filtration by closed subvarieties given by more restrictive rank conditions. We study the ideals defining these subvarieties, and products of such ideals, relating the geometry of rank loci with representation theory and commutative algebra. The ordinary powers of these ideals record products of equations, while their symbolic powers retain the primary components of the relevant associated primes. To understand the difference between the two constructions one is led to consider primary decomposition, minimal degrees of equations with prescribed symbolic vanishing, and containment between symbolic and ordinary powers. These questions have been studied mainly within regular rings, but determinantal rings provide a natural setting in which to study such questions inside singular ambient spaces.

These problems have a long history over polynomial rings. De Concini, Eisenbud, and Procesi developed a shape calculus for generic determinantal ideals and their products, including descriptions of primary decompositions, symbolic powers, and integral closures \cite{DEP}; see also Bruns--Vetter \cite{BrunsVetter}.
Abeasis and Abeasis--Del Fra developed the corresponding
invariant-ideal and tableau theories for symmetric minors and Pfaffians \cite{Abeasis,AbeasisDelFra}. Similar descriptions for catalecticant rank ideals and their products are given by Conca's straightening theory for Hankel matrices \cite{ConcaHankel}, and the subsequent treatment of extended Hankel matrices by Nam \cite{Nam}. All these theories share a common feature: membership in the relevant ideals is controlled by numerical functions on the shapes of standard monomials, denoted $\gamma_j$ below. For the generic, symmetric, and alternating settings, Henriques--Varbaro organized this common structure in their study of test, multiplier, and invariant ideals in  \cite{HenriquesVarbaro}.

The comparison of symbolic and ordinary powers has also yielded numerical invariants measuring their asymptotic behavior. While the Waldschmidt constant reflects the growth of the initial degree sequence of the symbolic powers,  resurgence, introduced by Bocci--Harbourne \cite{BocciHarbourne}, measures the sharpness of symbolic-power containments. Guardo--Harbourne--Van Tuyl have defined a similar invariant known as asymptotic resurgence in \cite{GuardoHarbourneVanTuyl2013}. This has been connected  to integral closures, Rees valuations, and skew Waldschmidt constants by DiPasquale--Francisco--Mermin--Schweig in \cite{DFMS}.
For classical rank ideals in polynomial rings, Kumar--Mukundan \cite{KumarMukundan} calculated initial degrees and resurgence invariants and analyzed containments of symbolic powers in products of ordinary powers with additional powers of the homogeneous maximal ideal. More recently, Bisui--Seceleanu
\cite{BisuiSeceleanu2026} determined symbolic powers and asymptotic invariants for arbitrary $\GL$-invariant ideals of a generic matrix and proved finite generation of their generalized symbolic Rees algebras.

The above-mentioned results for polynomial rings motivate us to study  the analogous problems in the coordinate ring of a bounded-rank variety. In this case, the ambient ring is a determinantal quotient itself, and products of rank prime ideals encode further rank conditions within matrices of bounded rank. Passing to this quotient raises two interrelated difficulties.
Imposing the rank bound given by the ambient ring eliminates standard monomials with large factors and can make components of the ambient primary decomposition redundant. This changes both the symbolic filtration and its asymptotic invariants. Our aim is to describe this change explicitly and uniformly in four settings:
\begin{enumerate}[label=\textup{(\roman*)}]
\item a generic matrix modulo its $(q+1)$-minors;
\item a symmetric matrix modulo its $(q+1)$-minors;
\item an alternating matrix modulo its $2(q+1)$-Pfaffians;
\item a Hankel matrix modulo its $(q+1)$-minors. 
\end{enumerate}

In each case let $A$ denote the quotient ring and
\[
 \fp_1\supset\fp_2\supset\cdots\supset\fp_q
\]
its chain of nonzero rank primes, with $\fp_1=A_+$.
We work over an algebraically closed field of characteristic zero.
For a Young diagram $\sigma=(\sigma_1\ge\cdots\ge\sigma_s\ge1)$ with
$\sigma_1\le q$, put
\[
 D_\sigma=\fp_{\sigma_1}\cdots\fp_{\sigma_s},
 \qquad
 \gamma_j(\sigma)=\sum_{i=1}^s(\sigma_i-j+1)_+,
 \qquad |\sigma|=\sum_{i=1}^s \sigma_i,
\]
where $(u)_+=\max\{u,0\}$.
Our results concern the products $D_\sigma$ of rank primes.

In the following theorem we summarize the main primary-decomposition, symbolic-power, and containment results of our paper. 
Once the formula in  \Cref{thm:intro-main}(i) has been established using the standard-monomial theory pertaining to each setting, the remaining parts of \Cref{thm:intro-main} then follow uniformly from the same combinatorial argument applied to all four rings.

\begin{theorem}\label{thm:intro-main}
With the notation above, set
\[
 \beta=\max\{1,|\sigma|-q(s-1)\},
 \qquad h=q-\beta+1,
 \qquad \fm=A_+.
\]
Then the following statements hold.
\begin{enumerate}[label=\textup{(\roman*)}]
\item There is a primary decomposition
\[
 D_\sigma=\bigcap_{j=1}^{\sigma_1}
                 \fp_j^{(\gamma_j(\sigma))}.
\]
In particular, $\fp_\beta$ is the unique inclusion-maximal
associated prime of $A/D_\sigma$.

\item For integers $m\ge1$, the symbolic powers with respect
to the full set of associated primes are
\[
 D_\sigma^{(m)}
       =\bigcap_{j=\beta}^{\sigma_1}
                    \fp_j^{(m\gamma_j(\sigma))}.
\]

\item For every integer $m\ge1$, the initial degree of the respective symbolic power is 
\[
 \alpha(D_\sigma^{(m)})
       =m\gamma_\beta(\sigma)
          +(\beta-1)
             \left\lceil\frac{m\gamma_\beta(\sigma)}h\right\rceil.
\]
For integers $m,r\ge1$ and $N\ge0$ one has
\[
 D_\sigma^{(m)}\subseteq\fm^N D_\sigma^r
 \quad\Longleftrightarrow\quad
 \alpha(D_\sigma^{(m)})\ge N+r|\sigma|.
\]

\item The Waldschmidt constant, resurgence, and asymptotic
resurgence respectively satisfy
\[
 \widehat\alpha(D_\sigma)
       =\frac{q\gamma_\beta(\sigma)}h,
 \qquad
 \rho(D_\sigma)=\widehat\rho(D_\sigma)
       =\frac{|\sigma|h}{q\gamma_\beta(\sigma)}
       =\frac{\alpha(D_\sigma)}{\widehat\alpha(D_\sigma)}.
\]

\item The symbolic Rees algebra
\[
 \cR_s(D_\sigma)
   =\bigoplus_{m\ge0}D_\sigma^{(m)}T^m
\]
is a finitely generated $A$-algebra.
In the generic, symmetric, and alternating settings, it is
a normal domain with rational singularities; in particular, it is Cohen--Macaulay.
\end{enumerate}
\end{theorem}

The finite-generation assertion in part~\textup{(v)} is constructive. We write symbolic powers as sums of ordinary products of rank ideals in \Cref{thm:finite-sum-unbounded}. The multiplicities of the factors and the symbolic exponents are encoded in a lattice-point monoid corresponding to a rational polyhedral cone, whose Hilbert basis yields finite sets of algebra generators. Examples show that generators of higher symbolic degree can be necessary. Normality of the three classical settings follows from a description  of the symbolic Rees ring as an intersection of valuation rings. A multiplicity-free action of an appropriate group  then yields rational singularities in each classical setting through the theory of affine spherical varieties.

The paper begins by recalling some background concerning symbolic powers and related invariants in \cref{sec:asymptotic-background} and establishing the common standard-monomial input. This includes results on symbolic-power descent in \cref{sec:package}. We compute the irredundant primary decompositions and give two explicit formulas for the symbolic powers of $D_\sigma$ in \cref{sec:symbolic-prod}. We further construct generators for their symbolic Rees algebras in \cref{sec:Rees}. Finally we calculate the exact initial degrees in \cref{sec:asymptotic}, which we use to resolve the containment problem, and to determine the resulting asymptotic invariants.

\begin{acknowledgement}
A.M. was partially supported by a Simons Dissertation Fellowship and NSF grant DMS-2349623.
K.M.S. was partially supported by NSF grants DMS-2236983,
DMS-2601940, and DMS-2044833. A.S. was partially supported by DMS-2401482. S.S. was partially supported by DMS-2342256.
\end{acknowledgement}

\begin{aidisclosure}
Generative AI was used to assist with literature searches.    
\end{aidisclosure}

\section{Symbolic powers and asymptotic invariants}\label{sec:asymptotic-background}

Throughout the paper, $\kfield$ is an algebraically closed field of
characteristic zero and all rings are standard graded $\kfield$-domains.
The characteristic-zero hypothesis is used for the classification of
invariant ideals.  Several of the straightening and primary
 decomposition statements hold in greater generality.

\subsection{Symbolic powers}

Let $I$ be a homogeneous ideal in a Noetherian domain $A$.  We use
symbolic powers with respect to the full set of associated primes:
\begin{equation}\label{eq:symbolic-definition}
 I^{(m)}=\bigcap_{\fp\in\Ass_A(A/I)}
 \bigl(I^mA_\fp\cap A\bigr).
\end{equation}
In \eqref{eq:symbolic-definition}, it suffices to intersect over
the inclusion-maximal elements of $\Ass_A(A/I)$.
In particular, if the associated primes form a chain with
maximal element $\fp$, which will be the case for our ideals of interest, then
$I^{(m)}=I^mA_\fp\cap A$.

For a nonzero homogeneous ideal $I$, let
\[
 \alpha(I)=\min\{d:I_d\ne0\}.
\]
The {\em Waldschmidt constant} of $I$ is
\[
 \widehat\alpha(I)=\lim_{m\to\infty}\frac{\alpha(I^{(m)})}{m}
 =\inf_{m\ge1}\frac{\alpha(I^{(m)})}{m}.
\]
The equality follows from subadditivity; see for example \cite[Lemma~3.3 and Definition~3.5]{DiPasqualeNguyenSeceleanu2023}.
The {\em resurgence} and {\em asymptotic resurgence}, introduced in this context
in~\cite{BocciHarbourne} and \cite{GuardoHarbourneVanTuyl2013}, respectively, are
\[
 \rho(I)=\sup\left\{\frac mr:I^{(m)}\not\subseteq I^r\right\}
\]
and
\begin{equation}\label{eq:asymptotic-resurgence}
 \widehat\rho(I)=
 \sup\left\{\frac mr:
 I^{(mt)}\not\subseteq I^{rt}\text{ for all }t\gg0\right\}.
\end{equation}
The asymptotic resurgence is controlled by the rational cones of the
ordinary and symbolic Rees algebras; see also~\cite{DFMS}.  This viewpoint will be used
rather than an explicit list of Rees valuations.

\subsection{Young diagrams and the $\gamma$-functions}

A \emph{Young diagram} is a finite weakly decreasing sequence
$\sigma=(\sigma_1,\ldots,\sigma_s)$ of positive integers. We can also view $\sigma$ as a sequence of rows of boxes of lengths $\sigma_1,\ldots,\sigma_s$. We recall that $|\sigma|=\sum_i\sigma_i$,  $(a)_+=\max\{0,a\}$, and for $j\ge1$, 
\begin{equation}\label{eq:gamma-def}
 \gamma_j(\sigma)=
 \sum_{i=1}^{s}(\sigma_i-j+1)_+,
\end{equation}
which counts the number of boxes to the right of column $j-1$.
We also set $\gamma_{q+1}(\sigma)=0$ when $\sigma_1\leq q$. 
The dual or conjugate diagram to $\sigma$ is $\sigma^\vee=n(\sigma)$ defined by 
$n_j(\sigma)=\#\{i:\sigma_i\ge j\}$. It is a diagram with $\sigma_1$ rows. The elementary identity
\begin{equation}\label{eq:gamma-difference}
 \gamma_j(\sigma)-\gamma_{j+1}(\sigma)=n_j(\sigma)
\end{equation}
will be used repeatedly.

\begin{example}
Given below is a Young diagram $\sigma=(4,3,2)$  and its dual
diagram $\sigma^\vee=(3,3,2,1)$. We have, for example, $\gamma_{1}(\sigma) = 9$ and $\gamma_{3}(\sigma) = 3$.
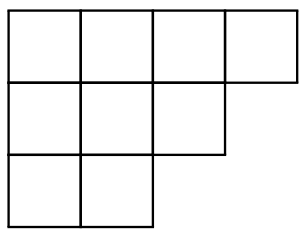
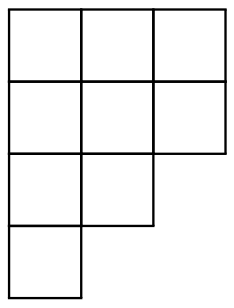
\[
\sigma=\ydiagram{4,3,2}
\qquad\qquad \qquad
\sigma^\vee=\ydiagram{3,3,2,1}
\]
\end{example}

If $\sigma$ and $\tau$ are diagrams, their concatenation
$\sigma*\tau$ is obtained by taking the multiset union of their parts and rearranging them in weakly decreasing order.  Thus
we write $\sigma^{*m}$ for the concatenation of $m$ copies of $\sigma$.

\section{The four determinantal settings}\label{sec:package}

This section records the representation theory and commutative algebra
needed later.  The emphasis is on the common shape calculus rather than
on the differences among the four straightening laws.

\subsection{Generic matrices}

Let $V$ and $W$ be vector spaces of dimensions $a$ and $b$,
respectively, and put $r=\min\{a,b\}$.  After choosing bases, let
$X=(x_{ij})$ be the generic $a\times b$ matrix, and
\[
 S_{\mathrm{gen}}=\kfield[X]
 =\Sym(V\otimes W^*).
\]
The group $G_{\mathrm{gen}}=\GL(V)\times\GL(W)$ acts naturally by $(M,N)\cdot X=MXN^{-1}$.  The  decomposition of  $S_{\mathrm{gen}}$ into irreducible $\kfield G_{\mathrm{gen}}$ modules and the invariant-ideal classification are as follows.
The decomposition written in terms of Schur functors  $\Schur$ applied to $V$ and $W^*$ respectively  is
\begin{equation}\label{eq:cauchy}
 S_{\mathrm{gen}}
 =\bigoplus_{\lambda:\,\lambda_1\le r}
 \Schur_{\lambda^\vee} V\otimes\Schur_{\lambda^\vee} W^*.
\end{equation}
We call the ideal $I_\lambda=( \Schur_{\lambda^\vee} V\otimes\Schur_{\lambda^\vee} W^*)$ of $S_{\mathrm{gen}}$ generated by the summand indexed
by $\lambda^\vee$ an {\em isotypic ideal}. For example if $\sigma=(t)$ is a diagram with one row then $\sigma^\vee=(1)^{*t}$ and  
\[
I_{(t)}=\left(\bigwedge^t V\otimes \bigwedge^t W^*\right)=I_t(X)
\] is the ideal of $t\times t$ minors of $X$.
By \cite[Theorem~4.1 and Proposition~5.1]{DEP} every homogeneous $G_{\mathrm{gen}}$-stable
ideal is a finite sum of  isotypic ideals. An important class of $G_{\mathrm{gen}}$-stable ideals consists of products of determinantal ideals. These are related to the isotypic ideals via integral closure. Indeed, 
\cite[Theorems~8.1]{DEP}  gives
\begin{equation}\label{eq:generic-isotypic-closure}
\overline{I_\lambda}
 =D_{\lambda}
 =I_{\lambda_1}(X)\cdots
 I_{\lambda_s}(X).
\end{equation}

We now discuss the analogous setup in determinantal rings. Fix $q<r$ and set
\begin{equation}\label{eq:generic-quotient}
 A_{\mathrm{gen}}=S_{\mathrm{gen}}/I_{q+1}(X),
 \qquad
 \fp_t=I_t(X)/I_{q+1}(X),\quad 1\le t\le q.
\end{equation}
The surviving summands in $A_{\mathrm{gen}}$ among those listed in  \eqref{eq:cauchy} are those indexed by
partitions with  $\lambda_1\leq q$.  We abuse notation and continue to write $I_\lambda$ for the isotypic ideal of $A_{\mathrm{gen}}$ generated by the corresponding summand and $D_\lambda$ for the image in $A_{\mathrm{gen}}$ of the eponymous ideal from \eqref{eq:generic-isotypic-closure}. 

Write a minor as $[i_1,\ldots,i_t\mid j_1,\ldots,j_t]$
with each index sequence strictly increasing and order them by
\[
 [I\mid J]\le[I'\mid J']
 \quad\Longleftrightarrow\quad
 t\ge t',\quad i_k\le i'_k,\ j_k\le j'_k\quad \text{ for } 1\le k\le t'.
\]
A standard bitableau is a product of minors $\delta_1\cdots\delta_s$ with
$\delta_1\le\cdots\le\delta_s$.  Equivalently, the two tableaux have strictly increasing rows
and weakly increasing columns.  Its shape is the Young diagram which lists the sizes of the
factors. These products form the standard basis
of $S_{\mathrm{gen}}$; thus modulo $I_{q+1}(X)$ the products with all sizes
at most $q$ form a basis of $A_{\mathrm{gen}}$. 

\subsection{Symmetric matrices}

Put $V=\kfield^{n}$, let $Y=(y_{ij})$ be the generic symmetric
$n\times n$ matrix ($y_{ij}=y_{ji}$ for all $1\leq i< j\leq n$), and 
\[
S_{\mathrm{sym}}=\kfield[Y]=\Sym(\Sym^2V).
\]
The natural congruence action of $\GL(V)$ given by $M\cdot Y=MYM^T$ gives the multiplicity-free
decomposition described by Abeasis~\cite[Theorem 2.4]{Abeasis} in terms of partitions whose columns are even, equivalently $\lambda^\vee\in 2\mathbb{N}^n$ as follows
\begin{equation}\label{eq:symmetric-decomp}
 S_{\mathrm{sym}}
 =\bigoplus_{\lambda:\,\lambda_1 \le n,\;\lambda^\vee\in 2\mathbb{N}^n}
 \Schur_{\lambda^\vee} V.
\end{equation}
The homogeneous $\GL(V)$-invariant ideals of $S_{\mathrm{sym}}$ are the sums of these irreducible
summands that are upward-closed under containment of Young diagrams; see \cite[Theorem 3.1]{Abeasis}.

Let $J_t(Y)$ be the ideal of $t$-minors of $Y$.  For $q<n$, set
\begin{equation}\label{eq:symmetric-quotient}
 A_{\mathrm{sym}}=S_{\mathrm{sym}}/J_{q+1}(Y),
 \qquad
 \fp_t=J_t(Y)/J_{q+1}(Y),\quad1\le t\le q.
\end{equation}
A Young diagram $\lambda$ with even column lengths has the form $\lambda=(\lambda_1 = \lambda_2 \geq  \lambda_3= \lambda_4 \geq  \lambda_5= \lambda_6, \ldots)$ and an even number $2s$ of parts. Thus the relevant information is preserved by
defining the following new diagram
\begin{equation}\label{eq:sym-transform}
 \delta_{\mathrm{sym}}(\lambda)
 =(\lambda_2,\lambda_4,
 \lambda_6,\ldots, \lambda_{2s}).
\end{equation}
We set  $D_{\delta_{\mathrm{sym}}(\lambda)}=J_{\lambda_2}(Y)J_{\lambda_4}(Y)\cdots J_{\lambda_{2s}}(Y)$ in $S_{\mathrm{sym}}$ or 
\[
D_{\delta_{\mathrm{sym}}(\lambda)}=\fp_{\lambda_2}\fp_{\lambda_4}\cdots \fp_{\lambda_{2s}}
\] in $A_{\mathrm{sym}}$. 
In particular, the partition $\lambda=(t, t)$ with two equal rows of length $t$ corresponds to $J_t(Y)$.

The use of $\SL(V)$ rather than $\GL(V)$ does not change the class of
homogeneous stable ideals.

\begin{proposition}\label{prop:sl-gl}
A homogeneous ideal of $A_{\mathrm{sym}}$ is $\SL(V)$-stable if and
only if it is $\GL(V)$-stable.
\end{proposition}

\begin{proof}
Only the forward implication requires proof.  Every element of
$\GL(V)$ is, after choosing an $n$-th root of its determinant, a product
of an element of $\SL(V)$ and a scalar matrix.  A scalar matrix acts by
a scalar on each homogeneous component of $S_{\mathrm{sym}}$ and hence
preserves every homogeneous ideal.  The same statement descends to the
homogeneous quotient $A_{\mathrm{sym}}$. Thus if an ideal of $A_{\mathrm{sym}}$ is both $\SL(V)$-stable and homogeneous, then it is $\GL(V)$-stable.
\end{proof}

The symmetric standard basis is described by a single tableau with
rows strictly increasing and columns weakly increasing satisfying the parity condition encoded in \eqref{eq:symmetric-decomp}. 
The same
functions \eqref{eq:gamma-def} govern symbolic powers, however in the $\gamma$-formulas we must use the normalized shape
$\delta_{\mathrm{sym}}(\lambda)$, listing one part for each minor,
not the doubled tableau shape; see \cite[Section~2.4]{HenriquesVarbaro}.

\subsection{Alternating matrices}

Put $V=\kfield^{n}$, let $Z=(z_{ij})$ be the generic alternating $n\times n$ matrix ($z_{ij}=-z_{ji}$ for all $1\leq i< j\leq n$ and $z_{ii}=0$ for all $1\leq i\leq n$), and 
\[
 S_{\mathrm{alt}}=\kfield[z_{ij}:1\le i<j\le n]
 =\Sym(\textstyle\bigwedge^2V).
\]
For $t\le\lfloor n/2\rfloor$, let $P_{2t}(Z)$ be the ideal generated by
all $2t$-Pfaffians of $Z$.  The natural congruence action of $\GL(V)$ given by $M\cdot Z=MZM^T$ gives rise to an  irreducible $\kfield\GL(V)$-decomposition of $S_{\mathrm{alt}}$, due to Abeasis--Del Fra~\cite[Theorem 2.5]{AbeasisDelFra}
\begin{equation}\label{eq:alternating-decomp}
 S_{\mathrm{alt}}
 =\bigoplus_{\lambda:\,\lambda_1 \le n,\;
                  \lambda_i\text{ even}}
 \Schur_{\lambda^\vee} V.
\end{equation}
Fix $q<\lfloor n/2\rfloor$ and set
\begin{equation}\label{eq:alternating-quotient}
 A_{\mathrm{alt}}=S_{\mathrm{alt}}/P_{2(q+1)}(Z),
 \qquad
 \fp_t=P_{2t}(Z)/P_{2(q+1)}(Z),\quad1\le t\le q.
\end{equation}
For a Young diagram $\lambda$ with even row lengths, put
\begin{equation}\label{eq:alt-transform}
 \delta_{\mathrm{alt}}(\lambda)
 =\left(\frac{\lambda_1}{2},
        \frac{\lambda_2}{2},\ldots, \frac{\lambda_s}{2} \right).
\end{equation}
Set
$
D_{\delta_{\mathrm{alt}}(\lambda)}
   =\prod_{i=1}^{s}P_{\lambda_i}(Z)
$
in $S_{\mathrm{alt}}$, and use the same notation for its image in $A_{\mathrm{alt}}$,
\[D_{\delta_{\mathrm{alt}}(\lambda)}=\prod_{i=1}^{s}\fp_{\lambda_i/2}.\]

A standard Pfaffian product has factors
$[i_1,\ldots,i_{2t}]$ ordered by decreasing $t$. Viewed as a tableau it has strictly increasing even size rows and weakly increasing columns. Its normalized shape lists the integers $t$, so its polynomial degree is the size of
that normalized shape.

Fixing a nondegenerate quadratic form $Q$ on $V$ gives an orthogonal
subgroup $\Orth(V,Q)\subset\GL(V)$.  Every Pfaffian-rank ideal and every
sum or product of such ideals is $\Orth(V,Q)$-stable.  The theorems of
this paper therefore apply to these $\Orth(V,Q)$-invariant ideals.  We
do not assert that every $\Orth(V,Q)$-stable ideal has this form.

\subsection{Hankel matrices}

Let $N\ge3$ and write $S_{\mathrm{Han}}=\kfield[x_1,\ldots,x_N]$.  For
$1\le t\le\lfloor(N+1)/2\rfloor$, consider the $t\times(N-t+1)$ Hankel
matrix
\begin{equation}\label{eq:hankel-matrix}
 \mathbf H_t=
 \begin{pmatrix}
 x_1&x_2&\cdots&x_{N-t+1}\\
 x_2&x_3&\cdots&x_{N-t+2}\\
 \vdots&\vdots&&\vdots\\
 x_t&x_{t+1}&\cdots&x_N
 \end{pmatrix}
\end{equation}
and let $K_t=I_t(\mathbf H_t)$ be its ideal of maximal minors. Every Hankel determinantal ring
is isomorphic to one where the defining ideal is generated by the maximal sized
minors of a Hankel matrix \cite[Lemma 2.3]{GrusonPeskine1982}; alternatively see \cite[Proposition 7]{Watanabe1997Hankel} or \cite[Corollary 2.2(b)]{ConcaHankel}. The
ideals $K_t$ form a chain
\[
 K_1\supset K_2\supset\cdots\supset
 K_{\lfloor(N+1)/2\rfloor}.
\]
For $q<\lfloor(N+1)/2\rfloor$, put
\begin{equation}\label{eq:hankel-quotient}
 A_{\mathrm{Han}}=S_{\mathrm{Han}}/K_{q+1},
 \qquad
 \fp_t=K_t/K_{q+1},\quad1\le t\le q.
\end{equation}
The variables may be viewed, up to  rescaling by appropriate  binomial coefficients, as the
coefficients of a binary form of degree $N-1$ 
\[
F(u,v)
=
\sum_{i=0}^{N-1}\binom{N-1}{i}x_{i+1} u^{N-1-i}v^i.
\]
The group $\GL_2(\kfield)$ acts naturally on binary forms by change of basis according to the formula
\[
(g\cdot F)(u,v)
=
F(au+bv,cu+dv)=
\sum_{i=0}^{N-1}\binom{N-1}{i}(g\cdot x_{i+1})
u^{N-1-i}v^i,
 \quad \text{for } g=
\begin{pmatrix}
a&b\\
c&d
\end{pmatrix}
\in \GL_2(\kfield).
\]
This  $\GL_2(\kfield)$-action preserves every $K_t$. Geometrically, the Hankel determinantal rings $A_{\mathrm{Han}}$ arise as homogeneous coordinate rings of higher order secant varieties of rational
normal curves, see \cite[Chapter 11.7]{Room}. In particular, the quotient
$S_{\mathrm{Han}}/K_{q+1}$ is the homogeneous coordinate ring of the $(q-1)$-st secant variety of the rational normal curve $C_{N-1}$ in $\mathbb{P}^{N-1}$, see \cite[Section 4]{Eisenbud1988}. 

Conca constructed a straightening law for the minors of Hankel
matrices~\cite[Theorem 2.6]{ConcaHankel}.  Standard products of Hankel minors have a
shape given by the minor sizes.  His theorem identifies the symbolic
powers of each $K_t$ and the primary decomposition of every power
$K_t^s$ by the same $\gamma$-functions as in the generic case.  The
full $\GL_2(\kfield)$-stable ideal theory involves plethysm and is not
multiplicity-free. 

\subsection{Symbolic-power descent}

The following formulas explain the descent of symbolic powers to the quotient rings $A$. In the generic determinantal case the formula below was noted in  \cite[Remark~10.5]{BrunsVetter}.

\begin{proposition}
\label{prop:classical-descent}
In the generic, symmetric, and alternating cases, let $P_j$
denote the ambient rank ideal, with $P_j=P_{2j}(Z)$ in the
alternating case, and put
$A=S/P_{q+1}$ and $\fp_j=P_j/P_{q+1}$.
For $1\le j\le q$ and $0\ne f\in A$, the local order
\[
 v_j(f)=\max\{a\ge0:f\in\fp_j^aA_{\fp_j}\},
\]
extends to a discrete valuation of $\operatorname{Frac}(A)$. We have
\begin{equation}\label{eq:quotient-symbolic-descent}
 \fp_j^{(m)}
 =\{f\in A:v_j(f)\ge m\}
 =\frac{P_j^{(m)}+P_{q+1}}{P_{q+1}}.
\end{equation}
\end{proposition}

\begin{proof}
Invert the upper-left $(j-1)$-minor, or the upper-left
$2(j-1)$-Pfaffian in the alternating case, using the empty
block when $j=1$. In the generic and symmetric cases, write
the matrix as
\[
 \begin{pmatrix}B&C\\ D&E\end{pmatrix},
 \qquad W=E-DB^{-1}C.
\]
Block elimination reduces the rank ideals to the corresponding
minor ideals of $W$; in the symmetric case $D=C^{\mathsf t}$
and $W$ is symmetric. Alternating block congruence gives the
analogous reduction of Pfaffian sizes by $2(j-1)$.
For the generic calculation see
\cite[Proposition~2.4]{BrunsVetter}.

After localizing at $\fp_j$, the eliminated parameters form
a coefficient field $F$, and $A_{\fp_j}$ becomes
$R_{\mathfrak n}$, where $R$ is the smaller standard graded
determinantal or Pfaffian domain of rank at most $q-j+1$,
and $\mathfrak n=R_+$.
Since
\[
\operatorname{gr}_{\fp_jA_{\fp_j}}A_{\fp_j}= \operatorname{gr}_{\mathfrak nR_{\mathfrak n}}
       (R_{\mathfrak n})\cong R,
\]

is a domain,
the local order is additive on nonzero products and satisfies
the valuation inequality on sums. It therefore extends by
$v_j(f/g)=v_j(f)-v_j(g)$.

The irreducible representation decompositions of the classical quotient rings $A$ are multiplicity-free;
see \cite[Sections~2.3--2.5]{HenriquesVarbaro} and the references
therein. Let $M_\lambda$ be the surviving irreducible summand
of normalized shape $\lambda$ in $A$, whose highest-weight vector is
\[
 K_\lambda=\delta_{\lambda_1}\cdots\delta_{\lambda_s}.
\]
Here $\delta_t$ is the upper-left $t$-minor, or the upper-left
$2t$-Pfaffian. If $t<j$, then $\delta_t$ is a unit at $\fp_j$.
If $j\le t\le q$, block reduction expresses it as a unit times
a nonzero homogeneous residual minor or Pfaffian of degree
$t-j+1$. Consequently,
\begin{equation}\label{eq:highest-weight-order}
 v_j(K_\lambda)
 =\sum_i(\lambda_i-j+1)_+
 =\gamma_j(\lambda).
\end{equation}
Since $\fp_j^{(m)}$ is group-stable, complete reducibility
and multiplicity-freeness show that it contains $M_\lambda$
exactly when it contains $K_\l$, that is, based on the above computation
if and only if 
 $\gamma_j(\lambda)\ge m$.

The same cutoff describes the ambient symbolic powers by
De Concini--Eisenbud--Procesi
\cite[Theorem~7.1]{DEP}, Abeasis
\cite[Theorem~5.1]{Abeasis}, and Abeasis--Del Fra
\cite[Theorem~5.1]{AbeasisDelFra}; see also
\cite[Theorems~2.1,~2.6, and~2.9]{HenriquesVarbaro}.
Passing to $A$ removes precisely the summands with
$\lambda_1>q$, proving \eqref{eq:quotient-symbolic-descent}.
The compatible standard bases give the basis assertion.
In the generic case, the descent identity is also stated
explicitly in \cite[Remark~10.5]{BrunsVetter}.
\end{proof}

For Hankel matrices, we use Conca's ambient straightening
theory \cite{ConcaHankel} and Nam's description \cite{Nam} of the diagonal initial ideals.
The following argument supplies the primarity needed to
pass from those ambient results to symbolic powers in the
rank quotient.

\begin{proposition}
\label{prop:hankel-descent}
Let $A=S/K_{q+1}$ be a Hankel rank quotient and
$\fp_j=K_j/K_{q+1}$. For $1\le j\le q$ and $m\ge1$,
\begin{equation}\label{eq:hankel-symbolic-descent}
 \fp_j^{(m)}
 =\frac{K_j^{(m)}+K_{q+1}}{K_{q+1}}.
\end{equation}
These ideals have as a basis the surviving standard products
of shapes $\lambda$ satisfying $\gamma_j(\lambda)\ge m$.
\end{proposition}

\begin{proof}
Put $Q=K_j^{(m)}+K_{q+1}$ and use the degree-lexicographic
order with $x_1>\cdots>x_N$. Conca's ambient straightening
and symbolic-power results~\cite{ConcaHankel}, in Nam's
initial-ideal formulation
\cite[Corollary~3.9, Theorem~3.12, and Corollary~3.21]{Nam},
give
\[
 J_t:=\operatorname{in}(K_t)
   =(x_{i_1}\cdots x_{i_t}:i_{u+1}\ge i_u+2),
 \qquad
 \operatorname{in}(K_j^{(m)})=J_j^{(m)}.
\]
The two summands of $Q$ are spanned by subsets of the same
standard basis, whose elements have pairwise distinct leading
monomials. Hence 
\begin{equation}\label{eq:hankel-initial-sum}
 \operatorname{in}(Q)=L:=J_j^{(m)}+J_{q+1}.
\end{equation}
We prove that $L$ is unmixed.

For $E\subseteq\{1,\ldots,N\}$, let $c(E)$ be the largest
size of a subset containing no consecutive indices, and set
$\Delta_h=\{E:c(E)\le h\}$.
The facets of $\Delta_h$ are unions of $h$ disjoint adjacent
pairs; in particular, they have $2h$ vertices and their
maximal consecutive runs have even lengths.
With $P_F=(x_i:i\notin F)$, we therefore have
\[
 J_j^{(m)}
 =\bigcap_{F\text{ facet of }\Delta_{j-1}}P_F^m.
\]
These facet and symbolic-power descriptions are given in
\cite[proof of Theorem~3.12]{Nam}.

Write $r=j-1$. We first observe that every face
$C\in\Delta_r$ contained in a facet $G$ of $\Delta_q$
extends to a facet $F$ of $\Delta_r$ contained in $G$.
Indeed, choose an inclusion-maximal
$C\subseteq C'\subseteq G$ with $c(C')\le r$.
If a consecutive run of $C'$ had odd length, it would have
an adjacent missing vertex inside $G$. Adjoining that vertex
would leave $c(C')$ unchanged, even if it joined two runs,
contradicting maximality. Thus all runs of $C'$ have even
length and $|C'|=2c(C')$.
If $c(C')<r$, a missing vertex could be added, increasing
$c$ by at most one. Hence $c(C')=r$, and $C'$ is the
required facet.

Every associated prime of the monomial ideal $L$ is generated
by variables and has the form $(L:M)$ for a monomial
$M\notin L$. Fix such an $M$.
Choose a facet $G$ of $\Delta_q$ containing
$E=\operatorname{supp}(M)$, and a facet $F_0$ of $\Delta_r$
such that
\[
 \sum_{i\notin F_0}\deg_{x_i}(M)<m.
\]
By the preceding observation, $F_0\cap E$ extends to a facet
$F\subseteq G$ of $\Delta_r$. Then
\[
 \sum_{i\notin F}\deg_{x_i}(M)<m.
\]
For every $i\in F$, the monomial $x_iM$ remains outside
$P_F^m$ and has support contained in $G$, so $x_iM\notin L$.
Consequently $(L:M)\subseteq P_F$.
Since $\sqrt L=J_j$ and $P_F$ is a minimal prime of $J_j$,
we obtain $(L:M)=P_F$.
Thus all associated primes of $S/L$ have dimension $2r$,
and $L$ is unmixed.

Here $\sqrt Q=K_j$. If $Q$ had an embedded associated
prime $\mathfrak r$, choose homogeneous $f$ with
$(Q:f)=\mathfrak r$. The modules
\[
 \frac{Q+(f)}Q
 \qquad\text{and}\qquad
 \frac{\operatorname{in}(Q+(f))}{L}
\]
have the same Hilbert function, since initial ideals preserve
Hilbert functions. The first has dimension
$\dim S/\mathfrak r<2r$, whereas the second is a nonzero
submodule of the unmixed module $S/L$ and therefore has
dimension $2r$. This contradiction proves that $Q$ is
$K_j$-primary.

Finally, the definition of the ambient symbolic power gives
\[
 (Q/K_{q+1})A_{\fp_j}
 =\fp_j^mA_{\fp_j}.
\]
Since $Q/K_{q+1}$ is $\fp_j$-primary, contraction proves
\eqref{eq:hankel-symbolic-descent}. The common standard
basis gives the asserted quotient basis.
\end{proof}

\subsection{The uniform standard-monomial package}

We now state the common input.  A \emph{standard monomial} means a
standard bitableau in the generic case, the corresponding symmetric or
Pfaffian tableau in the classical cases, and a standard product of
Hankel minors in the Hankel case.  Its shape is the diagram of its normalized 
factor sizes.

\begin{definition}\label{def:admissible-shape}
A diagram is \emph{admissible} if all its parts are at most $q$. 
\end{definition}

\begin{theorem}[Uniform standard-monomial package]\label{thm:package}
Let $A$ be one of
$A_{\mathrm{gen}},A_{\mathrm{sym}},A_{\mathrm{alt}},A_{\mathrm{Han}}$
and let $\fp_1\supset\cdots\supset\fp_q$ be the corresponding rank ideals described in the previous sections.
\begin{enumerate}[label=\textup{(\roman*)}]
 \item The ring $A$ is a normal Cohen--Macaulay domain, and each
 $\fp_j$ is prime.
 \item The standard monomials whose shapes have all parts at most $q$
 form a $\kfield$-basis of $A$.
 \item For $a\ge1$, the symbolic power $\fp_j^{(a)}$ has as a basis the
 standard monomials of shape $\lambda$ satisfying
 $\gamma_j(\lambda)\ge a$.
 \item If $\sigma$ is admissible, then
 \begin{equation}\label{eq:product-basis}
  D_\sigma=\prod_i\fp_{\sigma_i}
 \end{equation}
 has as a basis the standard monomials of shape $\lambda$ satisfying
 \begin{equation}\label{eq:all-gamma-ineq}
  \gamma_j(\lambda)\ge\gamma_j(\sigma)
  \quad\text{for }1\le j\le q.
 \end{equation}
 Equivalently,
 \begin{equation}\label{eq:uniform-primary}
  D_\sigma=\bigcap_{j=1}^{\sigma_1}
  \fp_j^{(\gamma_j(\sigma))}.
 \end{equation}
 \end{enumerate}
 \end{theorem}

\begin{proof}
Write $P_t$ for the ambient rank ideal, so that
$A=S/P_{q+1}$ and $\fp_t=P_t/P_{q+1}$.
Here $P_t=K_t$ in the Hankel case and
$P_t=P_{2t}(Z)$ in the alternating case.

For (i), the generic determinantal ring is a normal
Cohen--Macaulay domain by
\cite[Corollary~5.17 and Theorem~6.3]{BrunsVetter}.
The symmetric and alternating rings are respectively
orthogonal and symplectic invariant rings of polynomial
rings, by the first and second fundamental theorems
\cite{DeConciniProcesi1976}. They are normal domains,
and Cohen--Macaulayness follows from the
Hochster--Roberts theorem~\cite{HochsterRoberts}.
The Hankel assertion follows from
\cite[Theorem~2.1]{CMSV2018}.
Furthermore,
\[
 A/\fp_j\cong S/P_j
\]
is the corresponding lower-rank determinantal or Pfaffian
domain, or $\kfield$ when $j=1$. Thus every $\fp_j$ is prime.

The ambient standard bases are the classical straightening
bases~\cite{BrunsVetter,Abeasis,AbeasisDelFra} and
Conca's Hankel basis~\cite[Theorem~1.2]{ConcaHankel};
the latter's monomial indexing is also described in
\cite[Remark~3.18]{Nam}.
In each case, $P_{q+1}$ is spanned by precisely the standard
monomials having a factor of normalized size at least $q+1$.
Their removal gives the quotient basis in (ii).

Part (iii) is exactly the basis assertion of
\Cref{prop:classical-descent,prop:hankel-descent}.

For (iv), put
\[
 \mathcal D_\sigma=\prod_{i=1}^sP_{\sigma_i}\subseteq S.
\]
The ambient product and straightening theorems in the ring $S$ give a basis of $\mathcal D_\sigma$ consisting
of the standard monomials of shapes $\lambda$ satisfying
\[
 \gamma_j(\lambda)\ge\gamma_j(\sigma)
 \qquad\text{for every ambient rank index }j.
\]
For generic matrices this is
\cite[Theorem~10.9]{BrunsVetter}; the symmetric and
alternating results are due to
Abeasis~\cite{Abeasis} and
Abeasis--Del Fra~\cite{AbeasisDelFra}.
For Hankel matrices use
\cite[Theorem~3.12]{ConcaHankel}, or equivalently
\cite[Theorem~3.25]{Nam} together with
\cite[Corollary~3.21]{Nam}.

The image of $\mathcal D_\sigma$ in $A$ is $D_\sigma$.
Since $\mathcal D_\sigma$ and $P_{q+1}$ are spanned by
subsets of the same ambient standard basis, this image
has as a basis precisely the surviving standard monomials
satisfying the displayed inequalities. For indices $j>q$,
both $\gamma_j(\lambda)$ and $\gamma_j(\sigma)$ vanish,
so the surviving conditions are
\eqref{eq:all-gamma-ineq}.

By (iii), the intersection in \eqref{eq:uniform-primary}
is spanned by exactly these same basis elements.
This proves the equality. 
\end{proof}

\begin{remark}\label{rem:characteristic}
The generic statement is valid whenever the characteristic restrictions
in the primary-decomposition theorem are met; the Hankel straightening results are largely characteristic-free.  Characteristic zero is kept
throughout so that one may simultaneously use complete reducibility and
the invariant-ideal classifications as in the proof of \Cref{prop:classical-descent}.
\end{remark}

\section{Symbolic powers of products of rank prime ideals}\label{sec:symbolic-prod}

From this point on, the proof uses only
\Cref{thm:package}.  In particular, the four matrix types no longer
need to be considered separately. 

Throughout the section let $A$ be one of
$A_{\mathrm{gen}},A_{\mathrm{sym}},A_{\mathrm{alt}},A_{\mathrm{Han}}$
and let $\fp_1\supset\cdots\supset\fp_q$ be the respective  rank ideals described in the previous sections.
Moreover let $D_\sigma=\prod_{i=1}^s\fp_{\sigma_i}$ for an {\em admissible} diagram $\sigma$  with $s$ rows of length at most $q$.

\subsection{Primary decomposition after rank truncation}\label{sec:primary}

Fix an admissible diagram
$\sigma=(\sigma_1,\ldots,\sigma_s)$ with $\sigma_1\le q$.  
We first record an elementary comparison of partitions. We say that a partition  $\nu$ dominates $\mu$ if $\sum_{\ell=1}^{u} \nu_\ell\geq \sum_{\ell=1}^{u} \mu_
\ell$ for all integers $u\geq 0$.

\begin{lemma}\label{lem: hinge}
If $\mu,\nu$ are admissible partitions and $\nu$ dominates $\mu$ then for every $j$ such that $1\leq j\leq q$ we have $\gamma_j(\nu)\geq \gamma_j(\mu)$.
\end{lemma}
\begin{proof}
For a partition $\mu$, extended by zeros, and a real number $c\ge0$, we have
\begin{equation}\label{eq:hinge-sum}
 \sum_{\ell\ge1}(\mu_\ell-c)_+
 =
 \max_{u\ge0}
 \left\{\sum_{\ell=1}^{u}\mu_\ell-uc\right\}.
\end{equation}
Indeed, the maximum is attained when $u$ is the number of
parts of $\mu$ strictly larger than $c$. Consequently, if
$\nu$ dominates $\mu$, meaning that $\sum_{\ell=1}^{u} \nu_\ell\geq \sum_{\ell=1}^{u} \mu_\ell$ for all integers $u\geq 0$, then
\begin{equation}\label{eq:hinge-dominance}
 \sum_{\ell\ge1}(\nu_\ell-c)_+
 \ge
 \sum_{\ell\ge1}(\mu_\ell-c)_+
 \qquad\text{for every }c\ge0.
\end{equation}
Taking $c=j-1$ gives the desired conclusion.
\end{proof}

By \eqref{eq:uniform-primary}, $D_\sigma$ has a
possibly redundant primary decomposition indexed by $1\le j\le
\sigma_1$.  The next lemma determines precisely which components are needed.  It is analogous to \cite[Proposition 1.4]{DEP}.

\begin{lemma}[Irredundancy criterion]\label{lem:irredundancy}
The component $\fp_j^{(\gamma_j(\sigma))}$ in
\eqref{eq:uniform-primary} is irredundant if and only if
\begin{equation}\label{eq:irredundancy}
 \gamma_{j+1}(\sigma)
 \le(q-j)\bigl(n_j(\sigma)-1\bigr).
\end{equation}
\end{lemma}

\begin{proof}
Recall \eqref{eq:gamma-def} and \eqref{eq:gamma-difference} for the definitions of $\gamma_j$ and $n_j$.
If every part of $\lambda$ is at most $q$, then only the parts
satisfying $\lambda_i\ge j+1$ contribute to $\gamma_{j+1}(\lambda)$,
and each such part contributes at most $q-j$. Hence
\begin{equation}\label{eq:tail-bound}
\gamma_{j+1}(\lambda)
\leq
(q-j)n_{j+1}(\lambda)
\leq
(q-j)n_j(\lambda).
\end{equation}

Suppose that $\gamma_{j+1}(\sigma)>(q-j)(n_j(\sigma)-1)$. If a shape $\lambda$ satisfies
$\gamma_{j+1}(\lambda)\geq \gamma_{j+1}(\sigma)$,
then \eqref{eq:tail-bound} gives
\[
(q-j)n_j(\lambda)
\geq\gamma_{j+1}(\lambda)
\geq \gamma_{j+1}(\sigma)
>(q-j)(n_j(\sigma)-1),
\]
and therefore $n_j(\lambda)\geq n_j(\sigma)$.
Since $\gamma_j(\lambda)-\gamma_{j+1}(\lambda)=n_j(\lambda)$, we obtain
\[
\gamma_j(\lambda)
=
\gamma_{j+1}(\lambda)+n_j(\lambda)
\geq \gamma_{j+1}(\sigma) +n_j(\sigma)
=
\gamma_j(\sigma).
\]
Consequently, membership in  $\fp_{j+1}^{(\gamma_{j+1}(\sigma))}$ implies membership in $\fp_{j}^{(\gamma_{j}(\sigma))}$ under the  given assumption. 
Thus the $\fp_{j}$-primary component is redundant in \eqref{eq:uniform-primary}  under the stated
inequality.

Assume now that $\gamma_{j+1}(\sigma)\le(q-j)(n_j(\sigma)-1)$. We show that the $\fp_{j}$-primary component is irredundant in \eqref{eq:uniform-primary} by producing a standard monomial that belongs to every other primary component in \eqref{eq:uniform-primary}.

Put $n=n_j(\sigma)$ and $a=\gamma_{j+1}(\sigma)$.  Write the tail of the $n$ parts of
$\sigma$ which are at least $j$ as a partition of $a$ inside an
$n\times(q-j)$ rectangle.  Repack these $a$ boxes into an $(n-1)\times(q-j)$ rectangle so that
the resulting partition of $a$ is maximal in dominance order (equivalently,
fill the rows from left to right and top to bottom). The resulting tail
majorizes the old one and therefore does not decrease any
$\gamma_i$ with $i>j$ by \Cref{lem: hinge}.  Reattach $j$ boxes to the left of each of the $n-1$ rows of the new tail. Also reattach the original parts below row $n$ and add
two parts of size $j-1$ (when $j=1$ no such parts are needed) to obtain a new diagram $\lambda$.  As one box was discarded from row $n$ column $j$ of $\sigma$ and no new boxes were added to column $j$, the new
shape $\lambda$ satisfies
\[
 \gamma_j(\lambda)=a+n-1=\gamma_j(\sigma)-1,
\]
whereas $\gamma_i(\lambda)\ge\gamma_i(\sigma)$ for every $i\ne j$.
Indeed, for $i<j$ the loss from replacing $n$ parts crossing the
$j$-th column by $n-1$ such parts is $j-i+1$, while the two added
$(j-1)$-parts contribute $2(j-i)\ge j-i+1$.  

Thus the standard monomial of
shape $\lambda$ belongs to every primary component of \eqref{eq:uniform-primary} except the $\fp_j$-primary one.
Hence that component is irredundant.
\end{proof}

\begin{example}
The tail-repacking construction in the above proof is illustrated below for $\sigma=(5,4,3,2)$, $j=2$, $a=6$, $n=4$ and yields $\lambda=(5,5,2,1,1)$. The original and the repacked tails are shaded. Lost parts are displayed in red and gained parts in green.

\begin{figure}[h!]
    \centering
    \begin{tikzpicture}[>=stealth, baseline=(current bounding box.center)]

        \node(sigma) at (0,0) {%
            \begin{ytableau}
                \phantom{X} & \phantom{X} & *(gray)\phantom{X} &  *(gray)\phantom{X} &  *(gray)\phantom{X} \\
                \phantom{X} & \phantom{X} &  *(gray)\phantom{X} &*(gray) \phantom{X} \\
             *(red)   \phantom{X} & *(red)\phantom{X}&  *(gray)\phantom{X}  \\
                \phantom{X} & \phantom{X}
            \end{ytableau}
        };

        \node[above=0.25cm of sigma]
            {\(\text{diagram }\sigma\)};
        \node[below=0.2cm of sigma]
            {\(\displaystyle \sigma=(5,4,3,2)\subseteq 4\times5\)};
        
        \draw[->, very thick]
            (1.85,0) -- (3.6,0)
            node[midway,above=0.12cm]
            {\(\text{repack }\)};

        \node (repacked) at (5.5,0) {%
          \begin{ytableau}
                \phantom{X} & \phantom{X} & *(gray)\phantom{X} &  *(gray)\phantom{X} &  *(gray)\phantom{X} \\
                \phantom{X} & \phantom{X} &  *(gray)\phantom{X} &*(gray) \phantom{X} &  *(gray)\phantom{X} \\
                \phantom{X} & \phantom{X} 
            \end{ytableau}        };

        \node[above=0.25cm of repacked]
            {\(\text{repacked tail}\)};
        \node[below=0.2cm of repacked]
            {\(\displaystyle (5,5,2)\subseteq 3\times5\)};

        \draw[->, very thick]
            (7.35,0) -- (9.1,0)
            node[midway,above=0.12cm]
            {\(\text{add }(j-1)^2\)};

    
        \node (lambda) at (11,0) {%
           \begin{ytableau}
                \phantom{X} & \phantom{X} & *(gray)\phantom{X} &  *(gray)\phantom{X} &  *(gray)\phantom{X} \\
                \phantom{X} & \phantom{X} &  *(gray)\phantom{X} &*(gray) \phantom{X} &  *(gray)\phantom{X} \\
                \phantom{X} & \phantom{X} \\
              *(green)  \phantom{X}\\
             *(green)   \phantom{X}
            \end{ytableau}        };

\node[below=0.2cm of lambda]
            {\(\displaystyle \lambda=(5,5,2,1,1)\)};
           \node[above=0.25cm of lambda]
            {\(\text{new diagram}\)};

    \end{tikzpicture}

    \caption{The tail-repacking construction}
    \label{fig:tail-repacking}
\end{figure}
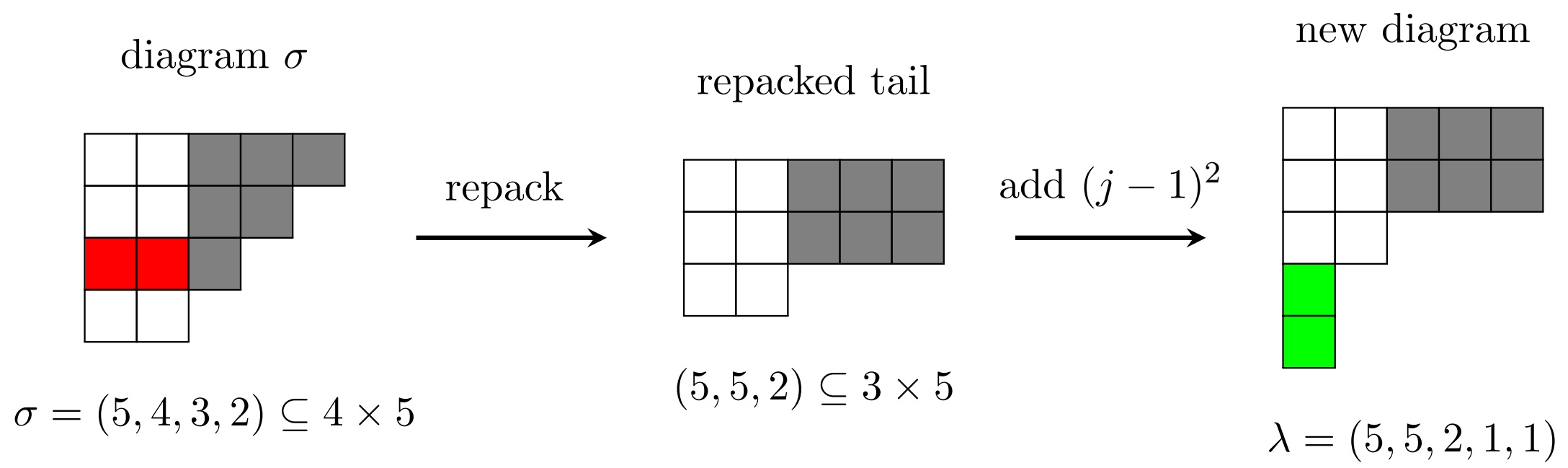
\end{example}

For the purpose of computing symbolic powers, we need to keep track of the maximal associated primes of ideals $D_\sigma$. Towards this end define
\begin{equation}\label{eq:beta-def}
 \beta_q(\sigma)=\max\{1,|\sigma|-q(s-1)\}.
\end{equation}
We usually write $\beta$ when $q$ and $\sigma$ are fixed.

\begin{proposition}\label{prop:max-ass}
The smallest index $j$ satisfying \eqref{eq:irredundancy} is
$\beta=\beta_q(\sigma)$.  Consequently,
\begin{align}
 \Ass_A(A/D_\sigma)
 &=\left\{\fp_j:1\le j\le\sigma_1,\;
 \gamma_{j+1}(\sigma)
 \le(q-j)(n_j(\sigma)-1)\right\},\label{eq:ass-set}\\
 \MaxAss_A(A/D_\sigma)&=\{\fp_\beta\}.\label{eq:maxass}
\end{align}
\end{proposition}

\begin{proof}
Suppose $\sigma$ is a partition with $s$ parts of length at most $q$. Let
\[
 d=\sum_{i=1}^{s}(q-\sigma_i),
 \qquad \text{ so that } \qquad
 |\sigma|-q(s-1)=q-d.
\]
If $q-d>1$, then $q-d\le\sigma_s$, so every part of $\sigma$ is at
least $j$ for $j\le q-d$.  In this range, $n_j(\sigma)=s$ and 
\[
 \gamma_{j+1}(\sigma)
 =s(q-j)-d.
\]
Thus \eqref{eq:irredundancy} is equivalent to $ s(q-j)-d \leq (q-j)(s-1)$, that is, 
$d\ge q-j$, or $j\ge q-d$.  Thus in this case the smallest index $j$ such that the $\fp_j$-primary component is irredundant in \eqref{eq:uniform-primary} is $\beta=q-d=|\sigma|-q(s-1)$.

If $q-d\le1$, the same calculation as above at
$j=1$ shows that the $\fp_1$-primary component is irredundant in \eqref{eq:uniform-primary}.  This proves that
$\beta=1$.  

The description of associated
primes follows from the irredundant primary decomposition, and the
chain $\fp_1\supset\cdots\supset\fp_q$ gives \eqref{eq:maxass}.
\end{proof}

\subsection{Formulas for symbolic powers of products of rank prime ideals}

\begin{theorem}\label{thm:symbolic-powers}
Let $\beta=\beta_q(\sigma)$ be as in \eqref{eq:beta-def}. For every integer $m\ge1$, the following identity holds
\begin{equation}\label{eq:symbolic-product-formula}
 D_\sigma^{(m)}
 =\bigcap_{j=\beta}^{\sigma_1}
 \fp_j^{(m\gamma_j(\sigma))}.
\end{equation}
Equivalently, $D_\sigma^{(m)}$ has a basis consisting of the standard
monomials of shapes $\lambda$ satisfying
\begin{equation}\label{eq:symbolic-shape-ineq}
 \gamma_j(\lambda)\ge m\gamma_j(\sigma)
 \quad\text{for }\beta\le j\le\sigma_1.
\end{equation}
\end{theorem}

\begin{proof}
By \Cref{prop:max-ass}, since the associated primes of $A/D_\sigma$ form a chain with maximal element $\fp_\beta$, the symbolic power may be computed by localizing
at $\fp_\beta$.  Since
$D_\sigma^m=D_{\sigma^{*m}}$ and $\gamma_j(\sigma^{*m})=m\gamma_j(\sigma)$, \Cref{thm:package}(iv) gives
\[
 D_\sigma^m
 =\bigcap_{j=1}^{\sigma_1}
 \fp_j^{(m\gamma_j(\sigma))}.
\]
For $j<\beta$, the prime $\fp_j$ strictly contains $\fp_\beta$, so its
component becomes the unit ideal after localization at $\fp_\beta$.
For $j\ge\beta$, the $\fp_j$-primary component  $\fp_j^{(m\gamma_j(\sigma))}$ is $\fp_\beta$-saturated, thus localization at $\fp_\beta$ and contraction to $A$ preserve the component.
Contracting back to $A$ gives \eqref{eq:symbolic-product-formula}.
The basis description is a consequence of \Cref{thm:package}(iii).
\end{proof}

In \Cref{thm:symbolic-powers} we gave an intersection formula for symbolic powers.  We now give a sum formula which, in addition,
produces explicit generators for the symbolic Rees algebra. This argument is independent of any
SAGBI or toric-degeneration assertions.

Set $d=q-\beta+1$. For $a=(a_\beta,\ldots,a_q)\in\mathbb{N}^d$, define
\[
  \tau(a)=\bigl(q^{a_q},(q-1)^{a_{q-1}},\ldots,\beta^{a_\beta}\bigr),
  \qquad
  D_{\tau(a)}=\prod_{t=\beta}^{q}\mathfrak{p}_t^{a_t},
  \qquad
  |a|=\sum_{t=\beta}^{q}a_t.
\]

\begin{theorem}\label{thm:finite-sum-unbounded}
Let $\beta=\beta_q(\sigma)$ be as in \eqref{eq:beta-def}. For every integer $m\ge1$, the following identity holds
\begin{align}
  D_\sigma^{(m)}
    &=\sum_{\substack{a\in\mathbb{N}^d\\ \gamma_j(\tau(a)) \geq m\gamma_j(\sigma) \; \forall \beta\leq j\leq \sigma_1}}
    D_{\tau(a)},   \label{eq:finite-sum-unbounded}
    \end{align}
    \end{theorem}

\begin{proof}
  By \Cref{thm:package}(iv),
\[
D_{\tau(a)}
  \subseteq \fp_j^{(\gamma_j(\tau(a)))}
  \subseteq \fp_j^{(m\gamma_j(\sigma))}
\]
whenever $\gamma_j(\tau(a))\ge m\gamma_j(\sigma)$. 
Intersecting over $j\in [\beta, \sigma_1]$ and applying  \Cref{thm:symbolic-powers} proves that every
summand in \eqref{eq:finite-sum-unbounded} is contained in the left side.

Conversely, let $\Delta$ be a standard product in the basis of $D_\sigma^{(m)}$ specified in  \Cref{thm:symbolic-powers}. For $\beta\leq t\leq q$, let
$a_t$ be the number of factors of normalized size $t$ in $\Delta$. Factors of size less than $\beta$
contribute zero to every $\gamma_j$, with $j\geq \beta$.  The product of the remaining factors
lies in $D_{\tau(a)}$, and
$\gamma_j(\tau(a))\ge m\gamma_j(\sigma)$ for
$\beta\le j\le\sigma_1$.
Consequently $\Delta\in D_{\tau(a)}$, proving
\eqref{eq:finite-sum-unbounded} by the basis assertion in  \Cref{thm:symbolic-powers}.
\end{proof}

\section{The symbolic Rees algebra of a product of rank prime ideals}\label{sec:Rees}

We now  discuss finite-generation, normality, and Cohen--Macaulayness of the symbolic Rees algebra for the ideals $D_\sigma$.

The proof uses convex geometry and is based on the following constructions. Fix a Young diagram $\sigma$ and integers $q>0$ and $\beta=\beta_q(\sigma)$. We introduce the integral matrix
\begin{equation}
  B=(b_{jt})_{\,\beta\leq j\leq \sigma_1,\,\beta\leq t\leq q},
  \qquad
  b_{jt}=(t-j+1)_+,
  \label{eq:B-matrix}
\end{equation}
and the vector 
\[
c=(\gamma_{\beta}(\sigma), \ldots, \gamma_{\sigma_1}(\sigma)).
\]
Consider the semigroup
\begin{equation}
  M_\sigma
  =\{(a,m)\in\mathbb{N}^d\times\mathbb{N}:Ba-mc\geq 0\}.
  \label{eq:monoid-M-sigma}
\end{equation}
It is the lattice point monoid of the pointed rational polyhedral cone
\[
  C_\sigma
  =\{(x,u)\in\mathbb{R}_{\geq 0}^{d+1}:Bx-uc\geq 0\}.
\] 
Therefore $M_\sigma$ is finitely generated by Gordan's lemma and saturated. Consider a (finite)  Hilbert basis for $M_\sigma$
\begin{equation}\label{eq: Hsigma}
  H_\sigma=\operatorname{Hilb}(M_\sigma),
  \qquad \text{ and put }
  H_\sigma^+=\{(a,m)\in H_\sigma:m>0\}.
\end{equation}

\begin{theorem}\label{thm:symbolic-rees}
Let $A$ be any of the generic, symmetric, alternating, or Hankel
rank-$q$ quotient rings considered above, and let $\sigma$ be an
admissible shape. Then for  the product $D_\sigma$ of rank ideals the symbolic Rees algebra satisfies
\begin{equation}
 \cR_s(D_\sigma)
 :=\bigoplus_{m\ge0}D_\sigma^{(m)}T^m =A\bigl[D_{\tau(a)}T^m:(a,m)\in H_\sigma^+\bigr].
\end{equation}
In the generic, symmetric, and alternating cases the symbolic Rees algebra is normal, has rational singularities, and hence is a normal Cohen--Macaulay domain. 
\end{theorem}

\begin{proof}
With notation as in \eqref{eq:B-matrix} we have  $(Ba)_j=\gamma_j(\tau(a))$. Consequently one can rewrite the set of inequalities $\gamma_j(\tau(a)) \geq m\gamma_j(\sigma)$ for $\beta_q(\sigma)\leq j\leq \sigma_1$ as $Ba\geq mc$. Here inequalities between vectors are understood
componentwise.

Then \eqref{eq:finite-sum-unbounded} gives 
\[
\mathcal{R}_s(D_\sigma)
  =A\bigl[
       D_{\tau(a)}T^m:(a,m)\in M_\sigma  \bigr].
\]
and the definition of $H_\sigma^+$ in \eqref{eq: Hsigma} further gives
\[
 A\bigl[
       D_{\tau(a)}T^m:(a,m)\in M_\sigma  \bigr] =A\bigl[
       D_{\tau(a)}T^m:(a,m)\in H_\sigma^+ \bigr]
\]
in view of the crucial identity
\[
D_{\tau(a+a')}=D_{\tau(a)}D_{\tau(a')}.
\]

In particular, taking minimal generators for each of the ideals $D_{\tau(a)}$ with $(a,m)\in H_\sigma^+$ gives a finite set of $A$-algebra generators for the symbolic Rees algebra.

For normality, let $K=\operatorname{Frac}(A)$ and put
$c_j=\gamma_j(\sigma)$. 
By \cite[Ch.~VI, \S10, no.~1, Lemma~1]{Bourbaki}, the valuation $v_j$ from \Cref{prop:classical-descent} extends to a valuation $w_j$ of
$K(T)$ with $w_j(T)=-c_j$ given by
the following formula on $A[T]$
\[
w_j\!\left(\sum_m f_mT^m\right)
   =\min_m\{v_j(f_m)-mc_j\}.
\]
The formula above extends to rational functions.
To see multiplicativity directly, choose a uniformizer $\pi$ for $v_j$, rescale $T$ by its $c_j$-th power, i.e., set $U=\pi^{c_j}T$ and observe that the coefficient of $U^m$ in the polynomial $\sum_m f_mT^m$ now has valuation $v_j(f_m)-mc_j$. Thus the weighted formula becomes the ordinary minimum of coefficient valuations in the new variable $U$ and multiplicativity follows by reducing coefficients to the domain $\kappa[U]$ where $\kappa$ is the residue field of the valuation ring of $v_j$.
 If $V_{w_j}$ is its valuation
ring, \Cref{thm:symbolic-powers} gives the equality inside $K(T)$
\begin{equation}\label{eq:rees-valuation-intersection}
\cR_s(D_\sigma)=A[T]\cap\bigcap_{j=\beta}^{\sigma_1}V_{w_j}.
\end{equation}
Indeed, for a polynomial the conditions are exactly
$v_j(f_m)\ge mc_j$ for every coefficient and every retained $j$. The rings on the right are integrally closed. Their intersection is therefore integrally closed. Since $D_\sigma\ne0$, the fraction field
of the left-hand side in \eqref{eq:rees-valuation-intersection} is $K(T)$: if $0\ne f\in D_\sigma$, then $T=(fT)/f$ lies in that fraction field. This proves normality.

For the generic, symmetric, and alternating cases let $G$ be $G_{\mathrm{gen}}$ or $\GL(V)$,
respectively. Each $D_\sigma^{(m)}$ is $G$-stable and is a direct sum of distinct irreducible summands of $A$. Let $\kfield^*$ act with weight $m$
on the factor $T^m$. Then $\cR_s(D_\sigma)$ is multiplicity-free as a $G\times\kfield^*$-module: a given irreducible of $G$ occurs at most once
in a fixed $T$-degree, and different $T$-degrees have different torus characters. The finite generation and normality already proved show
that its spectrum is a normal affine spherical variety, by \cite[Theorem~2.1.2]{Perrin}. Spherical varieties  have rational singularities and are Cohen--Macaulay
\cite[Corollary~2.3.4]{Perrin}. 
\end{proof}

The following examples illustrate the sum formula and the Hilbert-basis
construction in \Cref{thm:symbolic-rees}.  They also
show that higher symbolic degrees can be genuinely necessary.

\begin{example}\label{ex:32}
Set $q=3$ and $\sigma=(3,2)$, which leads to
 $\beta=2$, $\gamma_2(\sigma)=3$, and
$\gamma_3(\sigma)=1$.  Writing $a=(a,b)$, equation \eqref{eq:finite-sum-unbounded} becomes
\[
 D_{(3,2)}^{(m)}
 =\sum_{\substack{a,b\ge0\\a+2b\ge3m\\b\ge m}}
   \fp_2^a\fp_3^b.
\]
Removing nonminimal summands gives the closed formula
\[
 D_{(3,2)}^{(m)}
 =\sum_{j=0}^{\lfloor m/2\rfloor}
 D_{(3^{m+j},2^{m-2j})}.
\]
Consequently,
\begin{equation}\label{eq:32-rees}
 \cR_s(D_{(3,2)})
 =A[D_{(3,2)}T,D_{(3,3,3)}T^2]
 =A[\fp_2\fp_3T,\fp_3^3T^2].
\end{equation}
The $T$-degree two generator is necessary whenever the indicated products are
nonzero: it contains elements of polynomial degree $9$, whereas
$D_{(3,2)}^2$ starts in degree $10$.  In the alternating setting the
degree-one generators are products of a $4$-Pfaffian and a
$6$-Pfaffian, while the degree-two generators are products of three
$6$-Pfaffians.
\end{example}

\begin{example}\label{ex:rectangular-generators}
Consider the rectangular cases $(q,\sigma)=(4,(3,3))$ and $(q,\sigma)=(5,(4,4))$. 
Put $r=q-1$, so $r=3$ or $r=4$, and let
$U=\fp_{r-1}$, $V=\fp_r$, and $W=\fp_{r+1}$.  Here
$\beta=r-1$, and for multiplicities $(a,b,c)$ for these prime ideals the feasibility
conditions in \eqref{eq:monoid-M-sigma} are
\begin{equation}\label{eq:rect-constraints}
 a+2b+3c\ge4m,\qquad b+2c\ge2m.
\end{equation}
A Hilbert basis of the auxiliary monoid $M_\sigma$ in \eqref{eq:monoid-M-sigma} has positive-$m$ elements
\[
 (1,0,1;1),(0,2,0;1),(0,1,1;1),(0,0,2;1),
 (0,1,2;2),(0,0,3;2),(0,0,4;3).
\]
At the ideal level, nesting $W\subseteq V$ and the straightening
containment $UW\subseteq V^2$ remove several of the corresponding ideal products which would produce redundant generators for the symbolic Rees ring. Thus
\[
 \cR_s(D_{(r,r)})
 =A[V^2T,VW^2T^2,W^4T^3].
\]
Equivalently,
\[
 D_{(r,r)}^{(m)}
 =\sum_{u+2v+3w=m}
 D_{((r+1)^{2v+4w},r^{2u+v})}.
\]
In particular, depending on the value of $q$ we obtain
\[
\begin{aligned}
 q=4:\quad
 \cR_s(D_{(3,3)})
 &=A[D_{(3,3)}T,D_{(4,4,3)}T^2,D_{(4,4,4,4)}T^3],\\
 q=5:\quad
 \cR_s(D_{(4,4)})
 &=A[D_{(4,4)}T,D_{(5,5,4)}T^2,D_{(5,5,5,5)}T^3].
\end{aligned}
\]
The $T$-degree two and $T$-degree three generators cannot be omitted.  Indeed, in the generic matrix case
$VW^2$ starts in polynomial degree $3r+2<4r$, the initial degree of
$(V^2)^2$, while $W^4$ starts in degree $4r+4$, which is smaller than
the initial degrees $6r$ and $5r+2$ of the two degree-three products
available from the first two families.  The same formulas apply in the symmetric, alternating, and Hankel settings after replacing the rank ideals by the corresponding minors or Pfaffians; for Hankel this uses the admissible rectangular starting shape.
\end{example}

\section{Skew Waldschmidt constants and asymptotic resurgence}\label{sec:asymptotic}

For a nonzero ideal $J\subseteq A$ and $1\le j\le q$ define 
\[
\nu_j(J)=\max\{a\in\mathbb N:J\subseteq\fp_j^{(a)}\},
\qquad \fp_j^{(0)}=A.
\]
Equivalently, if $J$ has a standard basis and 
$1\le j\le q$, set
\[
 \nu_j(J)=\min\{\gamma_j(\lambda):
   \text{$J$ contains a standard monomial of shape $\lambda$}\}.
\]
The two definitions coincide since $\fp_j$-adic order is detected by the combinatorial function $\gamma_j$ by \Cref{thm:package}(iii).

Define the {\em skew Waldschmidt constants} associated to the functions $\nu_j$ as
\[
 \widehat\nu_j(D_\sigma)
 =\lim_{m\to\infty}\frac{\nu_j(D_\sigma^{(m)})}{m}.
\]
These limits will be computed explicitly below and therefore exist.
Since the polynomial degree of a standard monomial of shape $\lambda$
is $|\lambda|=\gamma_1(\lambda)$, one recovers the Waldschmidt constant as
\[
\widehat\nu_1(D_\sigma)=\widehat\alpha(D_\sigma).
\]

The standard monomial basis of $D_\sigma^r$ consists precisely
of the standard monomials whose shapes $\lambda$ satisfy
\[
\gamma_j(\lambda)\ge r\gamma_j(\sigma)
\qquad \text{ for }1\le j\le q.
\]
By contrast, \eqref{eq:symbolic-shape-ineq}  says that to find the standard monomial basis of $D_\sigma^{(m)}$  the indices $j<\beta_q(\sigma)$ are no longer explicit constraints. The next result computes how much of those missing inequalities is nevertheless forced by the surviving ones.

\begin{theorem}\label{thm:waldschmidt}
For an admissible diagram $\sigma=(\sigma_1, \ldots, \sigma_s)$ put $\beta=\beta_q(\sigma)$, $h=q-\beta+1$.
Then  for every $m\geq 1$ and  $1\le j\le q$ we have 
\begin{equation}\label{eq:exact-skew}
\nu_j(D_\sigma^{(m)})=\begin{cases}
m\gamma_\beta(\sigma)+(\beta-j)\left \lceil \frac{m\gamma_\beta(\sigma)}{h}\right \rceil & \text{ if } j<\beta\\
m\gamma_j(\sigma) & \text{ if } j \geq \beta.
\end{cases}
\end{equation}
In particular,
\begin{equation}\label{eq:exact-alpha}
\alpha(D_\sigma^{(m)})=m\gamma_\beta(\sigma)+(\beta-1)\left \lceil \frac{m\gamma_\beta(\sigma)}{h} \right \rceil.
\end{equation}
Consequently the
skew Waldschmidt constants of the product of rank ideals $D_\sigma$ are given by
\begin{equation}\label{eq:all-orders}
 \widehat\nu_j(D_\sigma)=
 \begin{cases}
 \displaystyle
 \frac{q-j+1}{q-\beta+1}\,\gamma_\beta(\sigma),&j<\beta,\\[8pt]
 \gamma_j(\sigma),&j\ge\beta.
 \end{cases}
\end{equation}
In particular,
\begin{equation}\label{eq:waldschmidt}
 \widehat\alpha(D_\sigma)
 =\frac{q\gamma_\beta(\sigma)}{q-\beta+1}.
\end{equation}
\end{theorem}

\begin{proof}
For $j\ge\beta$, \eqref{eq:symbolic-shape-ineq} gives
$\nu_j(D_\sigma^{(m)})\ge m\gamma_j(\sigma)$, while
$D_\sigma^m\subseteq D_\sigma^{(m)}$ and the shape $\sigma^{*m}$ in the standard monomial basis of  $D_\sigma^m$ gives
equality $\nu_j(D_\sigma^{(m)})= m\gamma_j(\sigma)$ and settles the last case of  \eqref{eq:exact-skew} and \eqref{eq:all-orders}.

Now suppose $j<\beta$ and let $\lambda$ be the shape of a standard monomial in $D_\sigma^{(m)}$. Put $E=\gamma_\beta(\l)$ and 
 $N=n_\beta(\lambda)=\#\{i:\lambda_i\geq \beta\}$. Each part counted
in $N$ contributes at most $h$ to $E$, and $E\ge m\gamma_\beta(\sigma)$ by  \eqref{eq:symbolic-shape-ineq}. Hence
\[
 N\ge\left\lceil E/h\right\rceil\ge\left\lceil m\gamma_\beta(\sigma)/h\right\rceil
 \qquad \text{ and }  \qquad
 \gamma_j(\lambda)\ge E+(\beta-j)N
 \ge m\gamma_\beta(\sigma)+(\beta-j)\left\lceil m\gamma_\beta(\sigma)/h\right\rceil.
\]

For equality, write $m\gamma_\beta(\sigma)=kh+r$, $0\le r<h$, and define a partition with $k$ or $k+1$ parts 
by
\begin{equation}\label{eq:extremal-shape}
 \theta=
 \begin{cases}
 (q^k),& \text{ if }r=0,\\
 (q^k,\beta-1+r),&  \text{ if } r>0.
 \end{cases}
\end{equation}
Deleting the first $\beta-1$ boxes of each row of $\sigma^{*m}$ gives
a partition of $m\gamma_\beta(\sigma)$ with parts at most $h$. The corresponding tail
of $\theta$ is $(h^k,r)$. This is
the dominance-maximal partition of $kh+r$ subject to the bound
$h$ on part size.  \Cref{lem: hinge} gives along the same lines as in the proof of \Cref{lem:irredundancy} that
\[
 \gamma_i(\theta)\ge m\gamma_i(\sigma)\quad \text{ for } i\ge\beta.
\]
Thus a standard monomial of shape $\theta$ belongs to
$D_\sigma^{(m)}$. It has $E=m\gamma_\beta(\sigma)$ and $N=\lceil m\gamma_\beta(\sigma)/h\rceil$, giving
simultaneous equality for every $j<\beta$. This proves
\eqref{eq:exact-skew}. Taking $j=1$ gives \eqref{eq:exact-alpha}. Dividing by $m$
and using $\lim_{m\to \infty}\lceil m\gamma_\beta(\sigma)/h\rceil/m= \gamma_\beta(\sigma)/h$ proves
\eqref{eq:all-orders}.
\end{proof}

Next we solve the containment problem and determine the asymptotic resurgence for products of prime ideals.

\begin{theorem}
Let $\sigma$ be an admissible partition. put
$\beta=\beta_q(\sigma)$ and $h=q-\beta+1$, and let
$\fm=A_+=\fp_1$. Then for all $m,r\ge1$ and $N\ge0$,
\begin{equation}\label{eq:exact-containment}
 D_\sigma^{(m)}\subseteq\fm^N D_\sigma^r
 \quad\Longleftrightarrow\quad
 m\gamma_\beta(\sigma)+(\beta-1)\left\lceil\frac{m\gamma_\beta(\sigma)}{h}\right\rceil\ge N+r|\sigma|.
\end{equation}
In particular, the following is a formula for the resurgence and asymptotic resurgence of a product of prime rank ideals
\begin{equation}\label{eq:resurgence-formula}
 \rho(D_\sigma)=\widehat\rho(D_\sigma)
 =\frac{|\sigma|h}{q\gamma_\beta(\sigma)}
 =\frac{\alpha(D_\sigma)}{\widehat\alpha(D_\sigma)}.
\end{equation}

\end{theorem}

\begin{proof}
Let $\nu=\sigma^{*r}*(1^N)$. Then $D_\nu=\fm^N D_\sigma^r$.
By \Cref{thm:package}(iv), a standard monomial of shape $\l$ is in $D_\nu$, if and only if it  satisfies the following inequalities
\begin{equation}\label{eq:target-inequalities}
 \gamma_1(\lambda)\ge N+r|\sigma|,\qquad
 \gamma_j(\lambda)\ge r\gamma_j(\sigma)\quad \text{ for } 2\le j\le q.
\end{equation}
Since the initial degree of the ideal $\fm^N D_\sigma^r$ is $\gamma_1(\nu)=N+r|\sigma|$ and the initial degree of $D_\sigma^{(m)}$ appears in \eqref{eq:exact-alpha},  the given containment implies the
inequality of initial degrees  in \eqref{eq:exact-containment}.

Conversely, assume that inequality and put $a_m=\lceil m\gamma_\beta(\sigma)/h\rceil$.
A nonzero product of shape $\sigma^{*m}$ shows
$\alpha(D_\sigma^{(m)})\le m|\sigma|$. The assumed degree inequality, together with the formula for \(\alpha(D_\sigma^{(m)})\) in \eqref{eq:exact-containment} gives 
$N+r|\sigma| \le \alpha(D_\sigma^{(m)})$. Combining these we deduce the inequality
\[
N+r|\sigma| \le m|\sigma|, \qquad \text{hence} \qquad r\le m.
\]

For $j\ge\beta$,
\Cref{thm:waldschmidt} gives
$\nu_j(D_\sigma^{(m)})=m\gamma_j(\sigma)\ge r\gamma_j(\sigma)$, where the inequality follows from that displayed above.
For $j<\beta$, use 
\[
\gamma_j(\sigma)=\gamma_\beta(\sigma)+(\beta-j)s\quad \text{ for }1\le j\le\beta
\]
 to obtain
\[
 \nu_j(D_\sigma^{(m)})-r\gamma_j(\sigma)
 =(m-r)\gamma_\beta(\sigma)+(\beta-j)(a_m-rs).
\]
The right side is an affine function of $x=\beta-j$ on
$[0,\beta-1]$. At $x=0$ it is nonnegative, and at $x=\beta-1$ it is
$\alpha(D_\sigma^{(m)})-r|\sigma|\ge N\ge0$. Hence it is nonnegative
throughout the interval. 
Together with the assumed degree inequality this proves that
\eqref{eq:target-inequalities} holds for standard monomials of shape $\l$ in $D_\sigma^{(m)}$. Since both $D_\sigma^{(m)}$ and  $\fm^N D_\sigma^r$ have standard monomial bases and we have shown that the standard monomials of the former are in the latter, containment follows.

Put $R=(|\sigma|h)/(q\gamma_\beta(\sigma))$ and $A_0=q\gamma_\beta(\sigma)/h$. If $m/r\ge R$, then
\[
\frac mr\ge R=\frac{|\sigma|}{A_0}
\quad\text{ so}\quad
mA_0\ge r|\sigma|,
\]
and
\[
\begin{aligned}
\alpha(D_\sigma^{(m)})
&=m\gamma_\beta(\sigma)
 +(\beta-1)\left\lceil\frac{m\gamma_\beta(\sigma)}h\right\rceil\\
&\ge m\gamma_\beta(\sigma)
 +(\beta-1)\frac{m\gamma_\beta(\sigma)}h\\
&=m\frac{q\gamma_\beta(\sigma)}h
=mA_0.
\end{aligned}
\]
Overall, $\alpha(D_\sigma^{(m)})\ge mA_0\ge r|\sigma|$, so by the previous criterion containment holds and
$\rho(D_\sigma)\le R$. 

If $m/r<R$, then
\[
\lim_{t\to \infty} \frac{\alpha(D_\sigma^{(mt)})}{t}= mA_0<r|\sigma|.
\]
Thus $D_\sigma^{(mt)}\not\subseteq D_\sigma^{rt}$ for all sufficiently
large $t$. Rational numbers $m/r<R$ approach $R$ from below, giving
$\widehat\rho(D_\sigma)\ge R$. Since
$\widehat\rho\le\rho$, \eqref{eq:resurgence-formula} follows after observing that
$\alpha(D_\sigma)=|\sigma|$.
\end{proof}

\end{document}